\documentclass{elsarticle}
\newtheorem{thm}{Theorem}[section]
\newtheorem{rem}{Remark}
\newtheorem{prop}{Proposition}[section]
\newtheorem{lemma}{Lemma}
\newtheorem{defn}{Definition}
\usepackage{bbold}
\usepackage{caption}
\usepackage{amsmath}
\usepackage{multirow}
\usepackage{amsfonts}
\usepackage{amssymb}
\usepackage{enumitem}
\usepackage{graphicx}
\usepackage{float}
\usepackage{rotating, color}
\usepackage{undertilde}
\usepackage{theoremref}
\allowdisplaybreaks
\numberwithin{equation}{section}

\DeclareMathOperator*{\argmax}{argmax} 

\usepackage{lineno,hyperref}
\modulolinenumbers[5]

\journal{Statistics and Probability Letters}

\biboptions{authoryear}

\begin{document}

\begin{frontmatter}

\title{Posterior Impropriety of some Sparse Bayesian Learning Models}


\author{Anand Dixit \corref{mycorrespondingauthor}}
\cortext[mycorrespondingauthor]{Corresponding author}
\ead{adixitstat@gmail.com}
\author{Vivekananda Roy}

\address{Department of Statistics, Iowa State University}

\begin{abstract}
Sparse Bayesian learning models are typically used for prediction in datasets with significantly greater number of covariates than observations. Such models often take a reproducing kernel Hilbert space (RKHS) approach to carry out the task of prediction and can be implemented using either proper or improper priors. In this article we show that a few sparse Bayesian learning models in the literature, when implemented using improper priors, lead to improper posteriors.  
\end{abstract}


\end{frontmatter}

\sloppy
\section{Introduction}

Modern datasets often have significantly greater number of covariates, $p$, than observations, $n$. For such datasets, often the objective is to predict the response variable for previously unobserved values of the covariates. If $p < n$, then one can fit a suitable linear model using a traditional statistical technique like ordinary least squares (OLS). But if $p > n$, then OLS is no longer applicable and hence one can rely on penalized methods such as least absolute shrinkage and selection operator (LASSO) proposed by \cite{tib::1996} or ridge regression proposed by \cite{hoe:ken::1970} to find a suitable model. But, both LASSO and ridge regression are penalized regression techniques that perform variable selection among the class of linear models. Hence, in case of $p > n$, if we wish to explore non linear class of models, we can estimate a function, $f$, from a functional space ($\mathcal{H}$) using the following Tikhonov regularization, 
\begin{equation}\label{eq:tik}
\min_{f \in \mathcal{H}} \Bigg[ \sum_{i=1}^n L(y_i, f(x_i)) + \lambda ||f||_{\mathcal{H}}^2 \Bigg],
\end{equation}
where $\{y_i, x_i\}_{i=1}^n$ is the training data such that $y_i \in \mathcal{R}$ for all $i$ and $x_i \in \mathcal{R}^p$ for all $i$, $L(\cdot, \cdot)$ is the loss function, $\lambda$ is the penalty parameter, $\mathcal{H}$ is the functional space and $||\cdot||_{\mathcal{H}}$ is the norm defined on $\mathcal{H}$. 

Since a functional space is infinite dimensional, the solution of \eqref{eq:tik} can also be infinite dimensional. Hence there is a possibility that we cannot use it for practical purposes. \cite{wah::1990} proved that, if the functional space is a reproducing kernel Hilbert space (RKHS), then the solution is finite dimensional and is given by, 
\begin{equation}\label{eq:tiksol}
f(x) =\sum_{j=1}^n k(x, x_j)  \beta_j,   
\end{equation}
where $k(\cdot, \cdot)$ is a reproducing kernel and $\{\beta_j\}_{j=1}^n$ are some unknown coefficients. The formal definition of RKHS and reproducing kernel can be found in \cite{ber:tho::2004}.  

\cite{tip::2001} used the finite dimensional solution \eqref{eq:tiksol} in a hierarchical Bayesian model to introduce the relevance vector machine (RVM). It has also been discussed in \cite{tip::2000} and \cite{bis:tip::2000}. The prior structure of RVM has been chosen in such a manner that it will produce a sparse solution and hence will lead to better predictions. RVM is a very popular sparse Bayesian learning model that is typically used for prediction and its popularity is demonstrated by large number of citations of the original RVM paper of \cite{tip::2001}.  

In Bayesian analysis, prior distributions are assumed on parameters. A prior distribution is said to be a proper prior if the prior density associated with it is a valid probability density function, else, it is said to be an improper prior. The most common objective prior used in the literature is the so-called Jeffreys' prior proposed by \cite{jef::1961} whose density function is directly proportional to the square root of the determinant of the Fisher information matrix, and, hence can be computed easily in several cases. The Jeffreys' prior can be proper or improper depending upon the data model used in the analysis. For Bayesian models involving improper priors, the posterior distribution of the parameters given the data is not guaranteed to be proper. Hence, in such cases, it is necessary to show that the normalizing constant associated with the posterior distribution is bounded above by a finite constant otherwise there is a possibility that the posterior distribution is improper and drawing inference from an improper posterior distribution is equivalent to drawing inference from a function that integrates to infinity.    
 
The RVM proposed by \cite{tip::2001} involves prior density functions proportional to $\lambda^{a-1}e^{-b\lambda}$ on different scale parameters. These densities can be assumed to be either proper gamma priors or improper priors based on the choice of the hyperparameters $a, b$ of that prior and \cite{tip::2001} presents both cases. The case of improper prior assumed by \cite{tip::2001} leads to an improper posterior distribution and we provide a proof of it. Additionally, we also derive necessary and sufficient conditions for the posterior propriety of RVM. The necessary conditions will help past researchers of RVM to check if the improper prior used by them leads to an improper posterior and the sufficient conditions will provide guidelines for future researchers to choose prior distributions that will guarantee posterior propriety. \cite{fig::2002} proposed to apply RVM using the popular Jeffreys' prior on parameters. The necessary conditions that we derive show that the choice of Jeffreys' prior also leads to an improper posterior. 

Sparse Bayesian learning models also involve classification models. \cite{mal::2005} proposed a RKHS based Bayesian classification model which makes use of the finite dimensional solution in \eqref{eq:tiksol} to build models corresponding to both logistic likelihoods as well as support vector machine related likelihoods. They propose to implement their model by using either proper priors or Jeffreys' prior. We show that the use of Jeffreys' prior in their models lead to an improper posterior. 

The article is structured as follows. In Section 2, we explain RVM and a related model proposed by \cite{fig::2002} along with their inference method in detail. Further in Section 2, we provide necessary and sufficient conditions for the posterior propriety of RVM and show that the sparse Bayesian learning models proposed by \cite{tip::2001} under improper prior and \cite{fig::2002} under Jeffreys' prior lead to improper posteriors. In Section 3, we provide details about the Bayesian classification models proposed by \cite{mal::2005} and show that the models are improper under the choice of Jeffreys' prior. Some concluding remarks are given in Section 4.              

\section{Relevance Vector Machine and its Impropriety}
Let $\{(y_i, x_i), i = 1, 2, \dots, n\}$ be the training data, where $y_i \in \mathcal{R}$ is the $i^{th}$ observation for the response variable and $x_i \in \mathcal{R}^p$ is the $p$ dimensional covariate vector associated with $y_i$. Let $y = (y_1, y_2, \dots, y_n)^T$ and $\beta=(\beta_0, \beta_1, \dots, \beta_n)^T$. Let $K$ be the $n \times (n+1)$ matrix whose $i^{th}$ row is given by $k_{i} = \big(1, k_{\theta}(x_i, x_1), k_{\theta}(x_i, x_2), \dots, k_{\theta}(x_i, x_n)\big)^T$ where $\{k_{\theta}(x_i, x_j): i=1, 2, \dots, n; j=1, 2, \dots, n\}$ are the values of the reproducing kernel and $\theta$ is a kernel parameter. The relevance vector machine proposed by \cite{tip::2001} is a hierarchical Bayesian model given as follows, 
\begin{subequations}
\label{eq:rvm}
\begin{align}
y|\beta, \sigma^2 &\sim N(K\beta, \sigma^2I),\label{eq:yrvm}\\
\beta|\lambda &\sim N(0, D^{-1})~~\text{with}~D = diag(\lambda_0, \lambda_1, \dots, \lambda_n),\label{eq:brvm}\\
\pi(\lambda_i) &\propto \lambda_i^{a - 1} \exp\{{-b\lambda_i}\} ~~for~all~i = 0, 1, 2, \dots, n, \label{eq:lrvm}\\
\pi\bigg(\frac{1}{\sigma^2}\bigg) &\propto \bigg(\frac{1}{\sigma^2}\bigg)^{c - 1} \exp\bigg\{-\frac{d}{\sigma^2}\bigg\}, \label{eq:srvm}
\end{align}
\end{subequations}
\noindent where ($a, b, c, d$) are user defined hyperparameters. Here $\{\sigma^2, \lambda_i: i = 0, 1, 2, \dots, n\}$ are assumed apriori independent. Also, $\beta$ and $\sigma^2$ are assumed apriori independent. The kernel parameter $\theta$ is typically estimated by cross validation in RVM. Let $\lambda = (\lambda_0, \lambda_1, \dots, \lambda_n)^T$. For $a > 0$ and $b > 0$, $\pi(\lambda_i)$ is a proper Gamma density with parameters $a$ and $b$ for all $i = 0, 1, 2, \dots, n$. Similarly, for $c > 0$ and $d > 0$, $\pi(1/\sigma^2)$ is a proper Gamma density with parameters $c$ and $d$. The posterior distribution of $(\beta, 1/\sigma^2, \lambda)$ corresponding to \eqref{eq:rvm} is given by,
\begin{equation}\label{eq:post}
\pi(\beta, 1/\sigma^2, \lambda|y) = \frac{f(y|\beta, \sigma^2) \pi(\beta, 1/\sigma^2, \lambda)}{m(y)},
\end{equation}  
\noindent where $f(y|\beta, \sigma^2)$ is the normal density in \eqref{eq:yrvm}, $\pi(\beta, 1/\sigma^2, \lambda)$ is the joint prior density of $(\beta, 1/\sigma^2, \lambda)$ derived from \eqref{eq:brvm}-\eqref{eq:srvm} and $m(y)$ is the marginal density defined as,
$$m(y) = \int_{\mathcal{R}^{n+1}_{+}}\int_{\mathcal{R}_{+}}\int_{\mathcal{R}^{n+1}} f(y|\beta, \sigma^2) \pi(\beta, 1/\sigma^2, \lambda) d\beta d\frac{1}{\sigma^2}d\lambda,$$

\noindent where $\mathcal{R}_{+} = (0, \infty)$. The posterior density given in \eqref{eq:post} is proper if and only if $m(y) < \infty$.  

The user defined hyperparameters can be chosen in such a way that the prior distribution imposed on the parameters turn out to be improper and in such cases the posterior propriety of the model is no longer guaranteed. The following theorems will provide necessary and sufficient conditions for the posterior propriety of RVM, that is, $m(y) < \infty$.

\begin{thm}\label{thm1}
Consider the RVM given in \eqref{eq:rvm}, then, for $b = 0$, which corresponds to prior $\pi(\lambda_i) \propto \lambda_i^{a - 1}~ for~all~i = 0, 1, \dots, n$, a necessary condition for the propriety of the posterior density \eqref{eq:post} is $a \in (-1/2, 0)$ for any choice of prior on $1/\sigma^2$, that is, for all $c, d \in \mathcal{R}$.    
\end{thm}
\noindent \textbf{Proof:} See the Appendix.
 
\begin{thm}\label{thm2}
Suppose $P_{K} = K(K^TK)^- K^T$ where $(K^TK)^-$ is a generalized inverse of $K^TK$. Then $(i)$ and $(ii)$ given below are the sufficient conditions for the propriety of the posterior density \eqref{eq:post}:   
\begin{enumerate}[label=(\roman*)]
\item The prior on $\lambda_i$ is a proper Gamma distribution for all $i = 0, 1, \dots, n, that~is,~a, b > 0$.
\item $y^T(I - P_{K}) y + 2d > 0$ and $c > -n/2$. 
\end{enumerate}
\end{thm}
\noindent\textbf{Proof:} See the Appendix.

\begin{rem}\label{rem1}
The proof of Theorem \ref{thm2} shows that posterior propriety of RVM is assured even if one wishes to use a proper prior on $\lambda_i$ other than Gamma distribution for all $i = 0, 1, \dots, n$ along with an improper prior on $1/\sigma^2$ that satisfies condition (ii) in Theorem \ref{thm2}.  
\end{rem}

Thus, Remark \ref{rem1} implies that posterior propriety is assured even if we choose the weakly informative half Cauchy prior on $\{\lambda_i^{-1/2}\}_{i=0}^n$ as suggested by \cite{gel::2006} or the type 2 Gumbel distribution (derived as the penalized complexity prior in \cite{sim::2017}) on $\{\lambda_i\}_{i=0}^n$, along with an improper prior like $\pi(1/\sigma^2) \propto \sigma^2$ or $\pi(1/\sigma^2) \propto 1$ for $1/\sigma^2$. Thus, Theorem \ref{thm2} and Remark \ref{rem1} provides researchers several other prior options for RVM that assures posterior propriety. 

In RVM, for given new values of the $p$ covariates, say, $x_{new}$, the objective is to predict the corresponding response variable, say, $y_{new}$. For predicting $y_{new}$, one can use the posterior predictive density given by, 
\begin{equation}\label{eq:predpost}
f(y_{new}|y) = \int_{\mathcal{R}^{n+1}_{+}}\int_{\mathcal{R}_{+}}\int_{\mathcal{R}^{n+1}}  f(y_{new}|\beta, \sigma^2)~\pi(\beta, 1/\sigma^2, \lambda|y) d\beta~d\frac{1}{\sigma^2}~d\lambda, 
\end{equation}
where $f(y_{new}|\beta, \sigma^2)$ is the density of $N(k^T_{new}\beta, \sigma^2)$ with $k_{new} = \big(1, k_\theta(x_{new}, x_1), k_\theta(x_{new}, x_2), \dots, k_\theta(x_{new}, x_n)\big)^T$ and $\pi(\beta, 1/\sigma^2, \lambda|y)$ is the posterior density defined in \eqref{eq:post}. \cite{tip::2001} approximated the posterior predictive density given in \eqref{eq:predpost} by,
$$\tilde{f}(y_{new}|y) = \int_{\mathcal{R}^{n+1}} f(y_{new}|\beta, \hat\sigma^2)~~\pi(\beta|\hat{\lambda}, \hat{\sigma}^2, y) d\beta,$$ 
where 
\begin{equation}\label{opt}
(\hat\lambda, \hat\sigma^2) = \argmax_{\lambda, \sigma^2} \pi(\lambda, 1/\sigma^2|y) = \argmax_{\lambda, \sigma^2} f(y|\lambda, \sigma^2),
\end{equation}
\noindent where $\pi(\lambda, 1/\sigma^2|y)$ is the marginal posterior density of $\lambda$ and $1/\sigma^2$, and, 
\begin{equation}\label{new}
f(y|\lambda, \sigma^2) = \int_{\mathcal{R}^{n+1}} f(y|\beta, \sigma^2)~~\pi(\beta|\lambda) d\beta.
\end{equation}
\noindent Using \eqref{eq:rvm}, simple calculations show that, 
\begin{align}
\beta|\hat\lambda, \hat\sigma^2, y~~&\sim~~N((K^TK + \hat{D}\hat{\sigma}^2)^{-1}K^Ty, (K^TK\hat{\sigma}^{-2} + \hat{D})^{-1})\nonumber\\
\implies~~y_{new}|y~~&\sim~~N(k_{new}^T(K^TK + \hat{D} \hat\sigma^2)^{-1}K^Ty, k_{new}^T(K^TK \hat{\sigma}^{-2} + \hat{D})^{-1}k_{new} + \hat{\sigma}^{2}). \nonumber 
\end{align} 
The mean of the above posterior predictive distribution is reported by \cite{tip::2001} as the predicted response when the observed covariates are $x_{new}$. In the above posterior predictive distribution used by \cite{tip::2001}, we also observe that the parameters $\lambda$ and $\sigma^2$ are estimated by maximizing the marginal density $f(y|\lambda, \sigma^2)$ and the prior imposed on them is $\pi(\lambda, \sigma^{-2}) \propto 1$ (Indeed the second equality in Equation \ref{opt} follows due to the use of this uniform prior.). Thus, the prior chosen is improper and is equivalent to choosing the hyperparameters $(a, b, c, d)$ in RVM, given in \eqref{eq:rvm}, to be $(1, 0, 1, 0)$. This choice of hyperparameters does not satisfy the necessary condition derived in Theorem \ref{thm1}. \cite{tip::2001} also mentions that optimizing $f(y|\lambda, \sigma^2)$ can be computationally challenging and hence he proposes to estimate $\log\lambda$ and $\log\sigma^{-2}$ by optimizing $\log f(y|\log\lambda, \log\sigma^{-2})$ and assuming uniform prior on $\log\lambda_i$'s and $\log\sigma^{-2}$, that is, $\pi(\log\lambda, \log\sigma^{-2}) \propto 1$, which is equivalent to $\pi(\lambda, \sigma^{-2}) \propto \sigma^2~\prod_{i=0}^n \lambda_i^{-1}$. Such a prior is also improper and can be formed by choosing the hyperparameters $(a, b, c, d)$ in \eqref{eq:rvm} to be $(0, 0, 0, 0)$. This choice of hyperparameters also violates the necessary conditions derived in Theorem \ref{thm1}. Thus, the RVM proposed by \cite{tip::2001} is based on an improper posterior. \cite{fig::2002} proposed to implement RVM by assuming the Jeffreys' prior on the prior variance parameters of $\beta$, that is, $\pi(1/\lambda_i) \propto \lambda_i$ for all $i$ which is equivalent to $\pi(\lambda_i) \propto 1/\lambda_i$ for all $i$. As mentioned before, this improper prior violates the necessary conditions derived in Theorem \ref{thm1}. Hence the model proposed by \cite{fig::2002} is  also based on an improper posterior. Thus, the necessary and sufficient conditions derived in Theorem \ref{thm1} and Theorem \ref{thm2} will be useful for past researchers to check if their choice of hyperparameters in RVM leads to a proper posterior. 

Interestingly, the prediction method of RVM can be viewed to be valid if $\{\lambda_i\}_{i=0}^n$ and $\sigma^2$ are assumed to be fixed at their estimates obtained by optimizing the marginal likelihood as given in \eqref{opt}. However, in the case of improper prior implementation of RVM, \cite{tip::2001} illustrates it as a Bayesian model in which flat improper priors are assumed on $\lambda$ and $1/\sigma^2$. Thus, there is a mismatch in the valid approach of prediction and model representation of RVM in the case of improper flat priors on $\lambda$ and $1/\sigma^2$. Using necessary and sufficient conditions derived in Theorem \ref{thm1} and Theorem \ref{thm2}, we hope to highlight this mismatch, and, provide greater clarity to RVM practitioners.     

\section{Sparse Bayesian Classification Model and its Impropriety}     
    
Let $y$ be an $n$ dimensional vector containing the observed response variables $\{y_i\}_{i=1}^n$ such that $y_i \in \{0, 1\}$ for all $i$ and let $z$ be an $n$ dimensional vector of latent variables that connect the response variables to the covariates. The Bayesian classification model based on reproducing kernels proposed by \cite{mal::2005} is as follows, 
\begin{align}
f(y|z) &\propto \exp\bigg\{-\sum_{i=1}^n l(y_i, z_i)\bigg\}\nonumber\\ 
z|\beta, \sigma^2, \theta &\sim N(K\beta, \sigma^2I)\nonumber\\
\beta|\lambda, \sigma^2 &\sim N(0, \sigma^2 D^{-1})~~\text{with}~D = diag(\lambda_0, \lambda_1, \dots, \lambda_n)\nonumber\\
\pi(\lambda_i) &\propto \lambda_i^{a - 1} \exp\{{-b\lambda_i}\}~for~all~i = 1, 2, \dots, n \nonumber\\
\sigma^2 &\sim IG(c, d)\nonumber\\
\theta&\sim U(u_1, u_2)\label{bclass}
\end{align}
where $y = (y_1, y_2, \dots, y_n)^T$, $z = (z_1, z_2, \dots, z_n)^T$, $l(\cdot, \cdot)$ is a loss function, $\beta=(\beta_0, \dots, \beta_n)^T$, $K$ is the $n \times (n+1)$ matrix whose $i^{th}$ row is given by $k_{i} = \big(1, k_{\theta}(x_i, x_1), k_{\theta}(x_i, x_2), \dots, k_{\theta}(x_i, x_n)\big)^T$ where $\{k_{\theta}(x_i, x_j): i=1, 2, \dots, n; j=1, 2,\dots, n\}$ are the values of the reproducing kernel, $\theta$ is the parameter in the reproducing kernel, $\lambda = (\lambda_0, \lambda_1, \dots, \lambda_n)^T$ with $\lambda_0$ fixed at a small number and ($a, b, c, d, u_1, u_2$) are user defined hyperparameters. For $X \sim IG(c, d)$, the density of the random variable $X$ is taken to be, $f(x) \propto~ x^{-c-1}e^{-d/x}~I(x >0)$ and $U(u_1, u_2)$ denotes the uniform distribution on the interval $(u_1, u_2)$. For $a>0$ and $b>0$, $\pi(\lambda_i)$ is a proper Gamma density with parameters $a$ and $b$. The parameters $\lambda_i$'s, $\sigma^2$ and $\theta$ are assumed apriori independent. In the case of Jeffreys' prior, the prior is assumed on $\lambda_0$ as well, that is, $\pi(\lambda) \propto \prod_{i=0}^n\lambda_i^{-1}$.  

The above model proposed by \cite{mal::2005} is quite general in nature, since it encompasses popular models like the logistic model and the support vector machine (SVM) model. \cite{mal::2005} recommend that the above model should be implemented using proper priors on $\lambda$ and $\sigma^2$ or by putting the Jeffreys' prior on $\lambda$ and a proper prior on $\sigma^2$. The following proposition shows that putting the Jeffreys' prior on $\lambda$ leads to an improper posterior.     
  
\begin{prop}\label{prop1}
If the Jeffreys' prior is assumed on $\lambda$ in the sparse Bayesian classification model given in \eqref{bclass}, that is, $\pi(\lambda) \propto \prod_{i=0}^n\lambda_i^{-1}$, then the posterior density of the parameters and latent variables of interest, $\pi(\beta, z, \sigma^2, \lambda, \theta|y)$ is improper.  
\end{prop}
\noindent\textbf{Proof:} See the Appendix.\\

\noindent Given new values of the $p$ covariates, say, $x_{new}$, sparse Bayesian classification model in \eqref{bclass} is used to predict the class $y_{new}$ belongs to. Since the response variable is binary, that is, $y_{new} \in \{0, 1\}$, the posterior predictive probability is given by, 
$$P(y_{new} = 1|y) = \int_{\Omega} P(y_{new} = 1|y, x_{new}, \omega) \pi(\omega|y)~d\omega$$
where $\Omega = \mathcal{R}^{n+1}\times\mathcal{R}^{n}\times\mathcal{R}_+\times\mathcal{R}^{n}_+\times(u_1, u_2)$, $\omega = (\beta, z, \sigma^2, \lambda, \theta)$, $\pi(\omega|y)$ is the posterior density of the parameters and latent variables of interest. Since the posterior distribution of the parameters and the latent variables of interest is not known in closed form, \cite{mal::2005} construct a Markov chain Monte Carlo (MCMC) sampler to draw samples from it and use those samples to produce a Monte Carlo estimate of $P(y_{new} = 1|y)$. If the Monte Carlo estimate is greater than 0.5, then $y_{new}$ is predicted to be 1 else 0. However, it is known that the usual Monte Carlo estimators converge to zero with probability one if the MCMC chain corresponds to an improper posterior distribution \citep{ath:roy::2014}. Also, generally MCMC samplers are incapable of providing a red flag when the posterior distribution is improper. In fact, \cite{hob:cas::1996} show that the MCMC draws from an improper posterior distribution may seem perfectly reasonable. Thus, to  detect posterior impropriety, one has to rely on theoretical analysis.  

\section{Conclusion and Discussion}
A probability density function is said to be valid only if the area under the curve is equal to one. This basic requirement is not assured for the posterior density function of a Bayesian model with an improper prior. Therefore for a Bayesian model with an improper prior, one should move ahead with inference only after showing that the posterior density function is valid. In this paper we have shown that some sparse Bayesian learning models based on improper priors do not have valid posterior density functions and therefore the inference or predictions drawn from them are not theoretically valid.  

In the case of hierarchical linear mixed models, \cite{rub:ste::2018} observed that posterior impropriety issues often arise due to assuming improper priors on parameters in the deeper levels of the hierarchy. Our observations concur with those of \cite{rub:ste::2018} as sparse Bayesian learning models considered in our paper assumed improper priors on parameters in the second level of the hierarchy which lead to improper posteriors. In the case of RVM, the sufficient conditions for posterior propriety derived in Theorem \ref{thm2} allows us to assume improper priors that satisfy condition $(ii)$ of Theorem \ref{thm2} on the parameter in the first hierarchical level as long as we assume a proper prior on parameters in the second hierarchical level. Thus, when assuming improper priors, it is important for Bayesian practitioners to establish posterior propriety using theoretical analysis. 

\section*{Acknowledgments}
The authors would like to thank the two referees and an Associate Editor for helpful comments that have improved the paper.  

\section*{Appendix}

\noindent\textbf{Notation.} Consider any two square matrices, say, $A$ and $B$ who have the same dimensions. Then, $A \le B$ means $B - A$ is a positive semidefinite matrix.   

\begin{defn}\label{def1}
Let $r = (r_1, r_2, \dots, r_n)^T\in\mathcal{R}^n$ and $s = (s_1, s_2, \dots, s_n)^T\in\mathcal{R}^n$ be any two $n$ dimensional vectors. A real valued function $f$ defined on $\mathcal{R}^n$ is said to be non decreasing in each of its arguments if $r << s$, that is, $r_i \le s_i~for~all~ i = 1, 2, \dots, n~\implies~f(r) \le f(s)$.    
\end{defn}

\begin{lemma}\label{lem2}
Let $P_{K} = K(K^TK)^-K^T$, where $(K^TK)^-$ is a generalized inverse of $K^TK$. 
Let $f_1(\lambda_0^{-1}, \lambda_1^{-1}, \dots, \lambda_n^{-1}) = \exp\bigg\{-\frac{1}{2}y^T\big(\sigma^2 I + KD^{-1}K^T\big)^{-1}y\bigg\}$. Then,\\ 
$$\exp\bigg\{-\frac{1}{2\sigma^2}y^Ty\bigg\} \le f_1(\lambda_0^{-1}, \lambda_1^{-1}, \dots, \lambda_n^{-1}) \le  \exp\bigg\{-\frac{1}{2\sigma^2}y^T(I - P_{K})y\bigg\}.$$

%
\end{lemma} 
\noindent\textbf{Proof:}\\
\noindent Differentiating $f_1$ with respect to $\lambda_i^{-1}$ we get, 
\begin{equation*}
\dfrac{\partial f_1}{\partial \lambda_i^{-1}} = \exp\bigg\{-\frac{1}{2}y^T\big(\sigma^2 I + KD^{-1}K^T\big)^{-1}y\bigg\}\frac{1}{2} y^T\big(\sigma^2 I + KD^{-1}K^T\big)^{-1}\big(KE_iK^T\big)\big(\sigma^2 I + KD^{-1}K^T\big)^{-1}y,
\end{equation*}
where $E_i$ is a $(n+1)\times (n+1)$ matrix with $1$ in the $i^{th}$ diagonal and $0$ everywhere else. Since $KE_iK^T$ is positive semidefinite, we get,\\  

\noindent $\dfrac{\partial f_1}{\partial \lambda_i^{-1}} \ge 0~~for~all~i~\implies f_1~\text{is a non decreasing function in each of its arguments.}$\\

\noindent Let $\lambda_{min} = min\{\lambda_0, \lambda_1, \dots, \lambda_n\}$ and $\lambda_{max} = max\{\lambda_0, \lambda_1, \dots, \lambda_n\}$. This implies that $\lambda^{-1}_{min} = max\{\lambda_0^{-1}, \lambda_1^{-1}, \dots, \lambda_n^{-1}\}$ and $\lambda^{-1}_{max} = min\{\lambda_0^{-1}, \lambda_1^{-1}, \dots, \lambda_n^{-1}\}$.  

\noindent Thus,
$$\lim_{\lambda^{-1}_{min}\rightarrow 0}~~f_1(\lambda_0^{-1}, \lambda_1^{-1}, \dots, \lambda_n^{-1}) = \exp\bigg\{-\frac{1}{2\sigma^2}y^Ty\bigg\}.$$
\noindent By Sherman - Morrison - Woodbury formula we have, 
\begin{equation}\label{eq:scomp}
\big(\sigma^2 I + KD^{-1}K^T\big)^{-1} = \frac{1}{\sigma^2}(I - K(K^TK + D\sigma^2)^{-1}K^T).
\end{equation}
\noindent Also,
\begin{equation}\label{eq:min}
(K^TK + D\sigma^2)^{-1} \leq (K^TK + \lambda_{min}\sigma^2~I)^{-1}.
\end{equation}
\noindent Using \eqref{eq:scomp} and \eqref{eq:min} we have, 
$$ f_1(\lambda_0^{-1}, \lambda_1^{-1}, \dots, \lambda_n^{-1}) \le \exp\bigg\{-\frac{1}{2\sigma^2}\bigg(y^T(I - K(K^TK + \lambda_{min}\sigma^2~I_{n+1})^{-1}K^T)y\bigg)\bigg\}.$$ 
\noindent By Lemma 1 of \cite{hob:cas::1996}, 
$$(K^TK)^{-} \equiv \lim_{\lambda^{-1}_{min}\rightarrow \infty}~~\bigg(K^TK + \frac{\sigma^2}{\lambda^{-1}_{min}}~I\bigg)^{-1}$$
\noindent is a generalized inverse of $K^TK$. Note that $\lambda^{-1}_{max}\rightarrow \infty$ implies $\lambda^{-1}_{min}\rightarrow \infty$. Hence we get,     
$$\lim_{\lambda^{-1}_{max}\rightarrow \infty}~~f_1(\lambda_0^{-1}, \lambda_1^{-1}, \dots, \lambda_n^{-1}) \le \exp\bigg\{-\frac{1}{2\sigma^2}y^T(I - P_{K})y\bigg\}.$$   
Thus, 
\begin{equation}\label{eq:exeq}
\exp\bigg\{-\frac{1}{2\sigma^2}y^Ty\bigg\} \le f_1(\lambda_0^{-1}, \lambda_1^{-1}, \dots, \lambda_n^{-1}) \le  \exp\bigg\{-\frac{1}{2\sigma^2}y^T(I - P_{K})y\bigg\}.
\end{equation}
\noindent The first inequality in \eqref{eq:exeq} also follows from the fact that $\sigma^2I + KD^{-1}K^T \ge \sigma^2I$. A similar argument could be used to prove the second inequality if $K^TK$ was non singular.

\begin{lemma}\label{lem4}
Consider the following integral,
\begin{equation}\label{eq:intg}
\int_{\mathcal{R}_+} \frac{t^{-(a+1)}}{(k + t)^{1/2}} dt
\end{equation}
where $k$ and $a$ are constants. The above integral is finite iff $a \in (-1/2, 0)$. In that case, the value of the integral is $ck^{-(a + 1/2)}$, where $c$ is some other constant.
\end{lemma} 
\noindent\textbf{Proof:}\\
Considering the transformation $t = k\tan^2\theta$, the integral in \eqref{eq:intg} becomes,
\begin{equation*}
2~k^{-(a + 1/2)}~\int_{0}^{\pi/2} (\sec^2\theta - 1)^{-(a+1)} \tan\theta\sec\theta~ d\theta.
\end{equation*}
Letting $z = \sec\theta$, the above integral becomes,
\begin{equation*}
2~k^{-(a + 1/2)}~\int_{1}^{\infty} (z^2 - 1)^{-(a+1)}~dz. 
\end{equation*}
The above integral is finite iff $a \in (-1/2, 0)$, thus proving the first part. Provided the above integral is some finite constant, say, $c/2$, the value of the integral given in \eqref{eq:intg} becomes $ck^{-(a + 1/2)}$. Thus proving the second part of the lemma.\\

{\large\noindent\textbf{Proof of Theorem \ref{thm1}}}\\ 

\noindent For RVM defined in \eqref{eq:rvm} with $b=0$,\\ 
\begin{footnotesize}
\begin{align}
f(y|\sigma^2) &= \int_{\mathcal{R}^{n+1}_+} \int_{\mathcal{R}^{n+1}} f(y|\beta, \sigma^2)~\pi(\beta|\lambda)~ \prod_{i=0}^n \lambda_i^{a - 1}~d\beta~d\lambda\nonumber\\
&=\int_{\mathcal{R}^{n+1}_+} \int_{\mathcal{R}^{n+1}} (2\pi)^{-n - 1/2}~\sigma^{-n}~|D|^{1/2} \exp\bigg\{-\frac{1}{2\sigma^2}\bigg((y - K\beta)^T(y - K\beta) + \beta^T D\sigma^2 \beta\bigg)\bigg\}~ \prod_{i=0}^n \lambda_i^{a - 1}~d\beta~ d\lambda\nonumber\\
&= \int_{\mathcal{R}^{n+1}_+} (2\pi)^{-n/2}~\sigma^{-n}~|D|^{1/2} |K^TK \sigma^{-2} + D|^{-1/2} \exp\bigg\{-\frac{1}{2}\bigg(\frac{y^Ty}{\sigma^2} - \frac{y^TK}{\sigma^2}(K^TK + D\sigma^2)^{-1}K^Ty\bigg)\bigg\}~ \prod_{i=0}^n \lambda_i^{a - 1}  ~d\lambda\nonumber\\
&= \int_{\mathcal{R}^{n+1}_+}\dfrac{\sigma}{(2\pi)^{n/2}}~|D|^{1/2}~|K^TK + D\sigma^2|^{-1/2}~\exp\bigg\{-\frac{1}{2\sigma^2}\bigg(y^T(I - K(K^TK + D\sigma^2)^{-1}K^T)y\bigg)\bigg\}~ \prod_{i=0}^n \lambda_i^{a - 1}~d\lambda\nonumber\\
&= \int_{\mathcal{R}^{n+1}_+}  \dfrac{\sigma}{(2\pi)^{n/2}}~|D|^{1/2}~|K^TK + D\sigma^2|^{-1/2}~\exp\bigg\{-\frac{1}{2}y^T\big(\sigma^2 I + KD^{-1}K^T\big)^{-1}y\bigg\} \prod_{i=0}^n \lambda_i^{a - 1} d\lambda,\nonumber\\
\end{align}
\end{footnotesize}

\noindent where $f(y|\beta, \sigma^2)$ is given in \eqref{eq:yrvm}, $\pi(\beta|\lambda)$ is the prior on $\beta$ given in \eqref{eq:brvm} and the last equality is obtained using \eqref{eq:scomp}. Let $e_1, e_2, \dots, e_{n+1}$ be the $n+1$ eigenvalues of $K^TK$ where $e_{max} = max\{e_1, e_2, \dots, e_{n+1}\}$. Then, $K^TK + D\sigma^2 \leq e_{max}~I + D\sigma^2$. Hence we get, 
\begin{equation}\label{eq:ineq1}
|K^TK + D\sigma^2|^{-1/2} \ge \prod_{i=0}^n \big(\lambda_i\sigma^2 + e_{max}\big)^{-1/2}.
\end{equation}
\noindent Using $|D|^{1/2} = \prod_{i=0}^n \lambda_i^{1/2}$, Lemma \ref{lem2}, \eqref{eq:ineq1} and letting $t = 1/\lambda_i$ for an arbitrary $i$, we get,   
$$f(y|\sigma^2) \ge  \dfrac{\sigma}{(2\pi)^{n/2}} \exp\bigg\{-\frac{1}{2\sigma^2}y^Ty\bigg\}\Bigg[\frac{1}{e_{max}^{1/2}}\int_{\mathcal{R}_+}\dfrac{t^{-(a+1)}}{\bigg(\frac{\sigma^2}{e_{max}} + t\bigg)^{1/2}}dt\Bigg]^{n+1}.$$
Using Lemma \ref{lem4}, the above integral is finite iff $a \in (-1/2, 0)$. Thus proving the necessary condition for the propriety of \eqref{eq:post} when $b=0$.\\

{\large\noindent\textbf{Proof of Theorem \ref{thm2}}}\\
 
\noindent For RVM defined in \eqref{eq:rvm}, 
\begin{align}
f(y|\sigma^2) &= \int_{\mathcal{R}^{n+1}_+} f(y|\lambda, \sigma^2) \pi(\lambda)d\lambda\nonumber\\
&=  \int_{\mathcal{R}^{n+1}_+}  \dfrac{\sigma}{(2\pi)^{n/2}}~|D|^{1/2}~|K^TK + D\sigma^2|^{-1/2}~\exp\bigg\{-\frac{1}{2}y^T\big(\sigma^2 I + KD^{-1}K^T\big)^{-1}y\bigg\} \pi(\lambda) d\lambda,\nonumber
\end{align}
where $\pi(\lambda)$ is the prior on $\lambda$ and $f(y|\lambda, \sigma^2)$ is given in \eqref{new}. 
Since $K^TK + D\sigma^2 \geq D\sigma^2$, we get, 
\begin{equation}\label{eq:ineq2}
|K^TK + D\sigma^2|^{-1/2} \le \prod_{i=0}^n \big(\lambda_i\sigma^2\big)^{-1/2}.
\end{equation}
\noindent Using $|D|^{1/2} = \prod_{i=0}^n \lambda_i^{1/2}$, Lemma \ref{lem2} and \eqref{eq:ineq2}, we get, 
\begin{equation}\label{eq:lamintg}
f(y|\sigma^2) \le \frac{1}{(2\pi)^{n/2}} \bigg(\frac{1}{\sigma^2}\bigg)^{n/2} \exp\bigg\{-\frac{1}{2\sigma^2}y^T(I - P_{K})y\bigg\} \int_{\mathcal{R}^{n+1}_+} \pi(\lambda)~ d\lambda.
\end{equation}
\noindent As mentioned before, as long as $\pi(\lambda)$ is a proper density, the integral in \eqref{eq:lamintg} will be 1. Therefore,   
\begin{equation*}
m(y) \le \frac{1}{(2\pi)^{n/2}}~ \int_{\mathcal{R}_+}\bigg(\frac{1}{\sigma^2}\bigg)^{n/2 + c - 1} \exp\bigg\{-\frac{1}{\sigma^2}\bigg(\dfrac{y^T(I - P_{K}) y}{2} + d\bigg)\bigg\}~d\frac{1}{\sigma^2}.  
\end{equation*} 
The integral above will be finite if $y^T(I - P_{K}) y + 2d  > 0$ and $c > -n/2$, thus proving the sufficient conditions for posterior propriety of RVM.\\

{\large\noindent\textbf{Proof of Proposition \ref{prop1}}}\\ 

\noindent For Bayesian classification model given in \eqref{bclass}, using similar calculations as in the proof of Theorem \ref{thm1}, we have,  
\begin{align}
f(z|\sigma^2, \theta) &= \int_{\mathcal{R}^{n+1}_+} \int_{\mathcal{R}^{n+1}}f(z|\beta, \lambda, \sigma^2, \theta)~\pi(\beta|\lambda) ~\pi(\lambda)~d\beta~d\lambda\nonumber\\
&=   \int_{\mathcal{R}^{n+1}_+}\dfrac{\sigma^{-n}}{(2\pi)^{n/2}}~|D|^{1/2}~|K^TK + D|^{-1/2}~\exp\bigg\{-\frac{1}{2\sigma^2}z^T\big(I + KD^{-1}K^T\big)^{-1}z\bigg\}\pi(\lambda)d\lambda.\nonumber
\end{align}
\noindent Since, $I + KD^{-1}K^T \geq I$, we get, 
\begin{equation}\label{eq:new}
 \exp\bigg\{-\frac{1}{2\sigma^2}z^T\big(I + KD^{-1}K^T\big)^{-1}z\bigg\} \geq \exp\bigg\{-\frac{1}{2\sigma^2}z^Tz\bigg\}.
\end{equation}
Further, $K^TK + D \leq e_{max}~I + D$ where $e_{max}$ is the largest eigenvalue of $K^TK$, hence we get, 
\begin{equation}\label{eq:ineq3}
|K^TK + D|^{-1/2} \ge \prod_{i=0}^n \big(\lambda_i + e_{max}\big)^{-1/2}.
\end{equation}
Using \eqref{eq:new}, \eqref{eq:ineq3} and letting $t = 1/\lambda_i$ for an arbitrary $i$, we get, 
$$f(z|\sigma^2, \theta) \ge \dfrac{\sigma^{-n}}{(2\pi)^{n/2}} \exp\bigg\{-\frac{1}{2\sigma^2}z^Tz\bigg\}\Bigg[\frac{1}{e_{max}^{1/2}}\int_{\mathcal{R}_+}\dfrac{t^{-1}}{\bigg(\frac{1}{e_{max}} + t\bigg)^{1/2}}dt\Bigg]^{n+1}.$$ 
From Lemma \ref{lem4}, the above integral is equal to $\infty$, thus proving Proposition \ref{prop1}.

\bibliography{Dixit_Roy_SPL}

\end{document}